\documentclass{amsart}
\usepackage{amssymb}
\usepackage{latexsym,calc,amssymb,amsfonts}
\usepackage{amsmath}
\usepackage{amscd}
\theoremstyle{plain}
\newtheorem{lemma}{Lemma}[section]
\newtheorem{thrm}{Theorem}[section]
\newtheorem{propn}{Proposition}[section]

\newtheorem{clm}{Claim}[section]

\theoremstyle{definition}

\newtheorem{define}{Definition}[section]

\newtheorem{remk}{Remark}[section]
\newtheorem{example}{Example}[section]
\newtheorem{algthrm}{Algorithrm}[section]

\newtheorem{conj}{Conjecture}

\DeclareMathOperator*{\aut}{GA} \DeclareMathOperator*{\aff}{Af}
\DeclareMathOperator*{\Gl}{GL} \DeclareMathOperator*{\El}{EA}
\DeclareMathOperator*{\BA}{BA} \DeclareMathOperator*{\tame}{T}

 \DeclareMathOperator{\tameequiv}{\ensuremath{\sim}}
 \DeclareMathOperator{\stameequiv}{\displaystyle \sim_{\bf{st}}}

\begin{document}
\title{Length Four Polynomial Automorphisms}
\author{Sooraj Kuttykrishnan}
\email{sooraj@cse.wustl.edu}
\address{Department of Computer Science,\ Washington University in St. Louis}

\begin{abstract}
We study the structure of length four polynomial automorphisms of $R[X,Y]$ when $R$ is a UFD. The results from this study are used to prove that if $\text{SL}_m(R[X_1,X_2,\ldots, X_n]) = \text{E}_m(R[X_1,X_2,\ldots, X_n])$ for all $n,\ m \ge 0$ then all length four polynomial automorphisms of $R[X,Y]$ that are conjugates are stably tame.
\end{abstract}
\maketitle
\section{Introduction}
Through out this paper $R$ will be a UFD. Amongst the many unanswered questions about the the structure of $\aut_2(R)$, the group of polynomial automorphisms of the polynomial algebra $R[X,Y]$, stable tameness conjecture is a long standing one. In this paper we will prove that certain length four automorphisms are stably tame. We will also give an intriguing example of a length four automorphism which is length four and stably tame.

 First we need a few definitions.
A polynomial map is a map $F=(F_1,..., F_n):
\mathbb{A}^n_R \rightarrow \mathbb{A}^n_R$ where each $F_i \in
R^{[n]}$. Such an $F$ is said to be invertible if there exists $G =
(G_1,..., G_n), G_i \in R^{[n]}$ such that $G_i(F_1,..., F_n)=X_i$
for $ 1\le i\le n$. The group of all polynomial automorphisms, $\aut_n(R)$ is defined as:
\begin{itemize}
\item $\aut_n(R)= \{F=(F_1,\ldots, F_n): F$ is invertible \}.
\end{itemize}
An important goal in the study of polynomial automorphisms is to understand the structure of this group in terms of some of its well understood subgroups.
An example of such a subgroup is
$$\text{Tame subgroup: }{\tame}_n(R)=\langle {\aff}_n(R), {\El}_n(R)\rangle \text{where}$$
\begin{align*}
{\aff}_n(R)=\{(a_{11}X_1+a_{12}X_2+\ldots +a_{1n}X_n+b_1,\ldots&,a_{n1}X_1+..a_{nn}X_n+b_n):\\
&(a_{ij}) \in {\Gl}_n(R) \:{\mbox{ and }}\:b_i \in R\}
\end{align*}
 is the subgroup of affine automorphisms of $\mathbb{A}^n_R$
and the elementary subgroup,
\begin{align*}{\El}_n(R)=\{\langle(X_1,X_2, \ldots , X_{i-1}, X_i+&f(X_1,\ldots,X_{i-1},\hat{X_i},
X_{i+1},\ldots , X_n),\ldots ,X_n)\rangle:\\
& f\in R[X_1,X_2,\ldots , \hat{X_i},\ldots, X_n],\ i \in \{1,\ldots, n\}\}
\end{align*}
Another well studied subgroup of $\aut_n(R)$ is the triangular subgroup,
\begin{align*}{\BA}_n(R)= \{\langle(a_1X_1+f_1(X_2,&\ldots , X_n), a_2X_2+f_2(X_3,\ldots, X_n), \ldots , a_nX_n+f_n)\rangle:\\
                                               &a_i \in R^*,\ f_i \in R[X_{i+1},\ldots ,X_n],\ 1\le i \le n-1, f_n \in R\}
 \end{align*}
 If $R$ is a domain, then $\aut_1(R)={\aff_1}(R)$. When $R$ is a field $k$ the following well known theorem gives us the structure of $\aut_2(k)$.
 \cite{Jung}, \cite{vanderKulk}
\begin{thrm}\label{Jungvanderkulk}(Jung, van der Kulk)
If k is a field then $\aut_2(k)= \tame_2(k)$. Further, $\tame_2(k)$
is the amalgamated free product of $\aff_2(k)$ and $\BA_2(k)$ over
their intersection.
\end{thrm}
A natural question that arises from Theorem \ref{Jungvanderkulk} is whether $\tame_3(k)$ is the whole group $\aut_3(k)$? Nagata
\cite{Nagata} conjectured that the answer is no and gave a candidate
counterexample.
\begin{example}\label{Nagatasexample}(Nagata)
\begin{align*}
\mbox{ Let } \:F_1&=(X,Y+{{\displaystyle{\displaystyle X^2\over t}}},t) \text{ and } F_2=(X+t^2Y,Y)\\
\mbox{ Then }\: N &=F_1^{-1}\circ F_2 \circ F_1.\\
                  &=(X+t(tY+X^2),Y-2(tY+X^2)X-t(tY+X^2)^2,t) \in \aut{_3}(k)
\end{align*}
\end{example}
Using the following algorithrm from \cite{Arnosbook}, we can conclude that $N\notin \tame_2(k[t])$.\\

Let $F=(P(X,Y),Q(X,Y))\in \aut_2(R)$ and $tdeg(F)=deg(P)+deg(Q)$ and $h_1$ be the highest degree term of $P$ and $h_2$ that of $Q$.
\begin{algthrm}\label{algorithrmtotesttameness}
Input: $F=(P,Q)$.\\
1) Let $(d_1,d_2)=(deg(P),deg(Q))$.\\
2) If $d_1=d_2=1$, go to 7.\\
3) If $d_1\neq d_2$, go to 5.\\
4) If there exists $\tau \in \aff_2(R)$ with $tdeg(\tau\circ F) < tdeg(F)$, replace $F$ by $\tau \circ F$ and go \indent to 1, else stop : $\notin \tame_2(R)$.\\
5) If $d_2<d_1$, replace $F$ by $(Q,P)$.\\
6) If $d_1\mid d_2$ and there exists $c\in R$ with $h_2= ch_1^{d_2/d_1}$, replace $F$ by $(X,Y-cX^{d_2/d_1})\circ F$ \indent and go to 1, else stop : $F \notin \tame_2(R)$.\\
7) If $\det JF \in R^*$, stop: $F\in \tame_2(R)$, else stop : $F\notin \tame_2(R)$.
\end{algthrm}
Shestakov and Umirbaev in 2002 \cite{NagataisWild} proved that $N
\notin \tame_3(k)$ and thus proved Nagata's conjecture.
\begin{define}
\

Let $ F,G \in \aut_n(R)$. Then
\begin{enumerate}
              \item $F$ is \emph{stably tame} if there exists $m \in {\mathbb
{N}}$ and new variables $ X_{n+1},\ldots, X_{n+m}$ such that the
extended map $\widetilde{F}=(F, X_{n+1},\ldots, X_{n+m})$ is tame.\\
 i.e $(F, X_{n+1},\ldots, X_{n+m})\in \tame_{n+m}(R)$

 \item  $F$ is \emph{tamely equivalent}${(\tameequiv)}$ to $G$ if there exists $H_1, H_2 \in
\tame_n(R)$ such that $ H_1\circ F \circ H_2 = G$.

\item $F$ is \emph{stable tamely equivalent}$(\stameequiv)$ to $H
\in \aut_{n+m}(R)$ if there exists $ \widetilde{H_1}, \widetilde{H_2
}\in \tame_{n+m}(R)$ such that $\widetilde{H_1} \circ \widetilde{F}
\circ \widetilde{H_2}=H$ where $\widetilde{F}=(F, X_{n+1},\ldots ,
X_{n+m})$
 \end{enumerate}
\end{define}
Martha Smith proved \cite{MarthaSmith} that $N$ from Nagata's example is stably tame with one more variable. This result led to the formulation of the following conjecture.
\begin{conj} If k is a field and $F \in \aut_n(k)$ then $F$ is stably tame.
\end{conj}
In her proof of the stable tameness of Nagata's example, Martha Smith exploited the decomposition of $N$ in Example \ref{Nagatasexample} into certain special type of elementary automorphisms as shown in the example. This led to further study of such decompositions and the notion of the length of an automorphism, which we discuss below.
The following proposition due to Wright is well known and a proof is given in \cite{Kuttykrishnanpaper1}.
\begin{propn}
Let $R$ be a domain $K$ its fraction field and $F\in \aut_2(R)$.
Then $F=L\circ D_{a,1}\circ F_{m}\circ F_{m-1}\circ...\circ F_1$
where $L=(X+c, Y+d),\ D_{a,1}=(aX,Y),\ F_i=(X,Y+f(X))\text{ or }
F_i=(X+g(Y),Y)$ for some $c,d \in R,\ a \in R^*,\ f(X),\ g(X) \in
K[X] $
\end{propn}
\begin{define}\
\begin{enumerate}
\item \emph{Length} of $F \in \aut_2^0(R)$ is the smallest natural number m such that
 $F=D_{a,1}\circ F_m\circ F_{m-1} \circ \ldots \circ F_2 \circ F_1$ where each $F_i$ is either
of the type $(X, Y+f_i(X))$ or $ (X+g_i(Y),Y)$ with $f_i(X),\ g_i(X)
\in K[X], a\in R^*$ and $f_i(0)=g_i(0)=0$.
\item $\text{L}^{(m)}(R)=\{F \in {\aut}_2^0(R):\ F \text { is of length } m\}$
\end{enumerate}
\end{define}
\begin{remk}\label{Da1doesn't matter}If $F \in \text{L}^{(m)}(R)$ as above and $F=D_{a,1}\circ F_m\circ F_{m-1} \circ \ldots \circ F_2 \circ F_1 \in \text{L}^{(m)}(R)$ then $F$ is tamely equivalent to $G=F_m\circ F_{m-1} \circ \ldots \circ F_2 \circ F_1$. Thus $F$ is stably tame iff $G$ is stably tame.
\end{remk}
It is easy to see that Nagata's example is of length three and it is stably tame with one
more variable. Drensky and Yu \cite{TameAndWild} began a systematic
study of length three automorphisms and proved the following result.

\begin{thrm}(Drensky, Yu) Let k be a field of characteristic zero and $ F \in \text{ L }^{(3)}(k[t])$ such that
 $F = F_1^{-1}\circ G \circ F_1$ where $F_1=(X,Y+f(X)), G=(X+g(Y),Y)$ with $f(X),\ g(X) \in k[t][X].$
  Then $F$ is stably tame with one more variable.
\end{thrm}
The following theorem was proved in \cite{Kuttykrishnanpaper1}.
Let $\text{SL}_n(R)$ denote the set of all $n\times n$ matrices with entries from R and determinant equal to 1 and $\text{E}_n(R)$ denote the group generated by the set of all nxn elementary matrices with entries from $R$.

\begin{thrm}\label{length3isstablytame1}Suppose $R$ is a UFD such that $\text{SL}_m(R[X_1,X_2,\ldots, X_n]) =\\ \text{E}_m(R[X_1,X_2,\ldots, X_n])$ for all $n,\ m \ge 0$. Then $F\in \text{L}^{(3)}(R) \Rightarrow F$ is stably tame.
\end{thrm}
This theorem was also claimed by Edo in \cite{EDOtotallystablytame} without the assumption that $\text{SL}_m(R[X_1,X_2,\ldots, X_n]) =\\ \text{E}_m(R[X_1,X_2,\ldots, X_n])$ for all $n,\ m \ge 0$. However, it is the author's contention that this assumption is required for the proof provided in \cite{EDOtotallystablytame} to hold.
So a natural question at this point to ask is if $F \in \text{L}^{(4)}(R)$ stably tame? As an evidence to an affirmative answer to this question, we prove the following theorem in this paper.
\begin{thrm}(Main Theorem)\label{length4commutatorisstablytame}\ Let R be a UFD and $F \in \text{L}^{(4)}(R)$ and $F = G_1^{-1}\circ F_1^{-1}\circ G_1 \circ
F_1$ where $ F_1=(X, Y+ f(X)),\ G_1= (X+g(Y),Y),\ f(X),\ g(X) \in
K[X] \mbox{ with } f(0) = g(0)= 0$. Then $F$ is stably tame.
\end{thrm}
\begin{remk}In \cite{berson-2007} Berson,van den Essen and Wright recently proved that if $F \in \aut_2(R)$, where $R$ is a regular ring then $F$ is stably tame. This is a much stronger result. However, our result does not require the ring to be regular.
\end{remk}
Before we present the proof of the main theorem, here is an example of a non tame automorphism of length four.
\begin{example}
Let $R$ be a UFD and $t\in R\backslash \{0\}$.
\begin{align*}
\text{Let }F_1&=(X,Y+\frac{(t+1)^3X^2}{t})\\
G_1&=(X+\frac{t^2Y}{(t+1)}) \text{ and }\\
F&=G_1^{-1}\circ F_1^{-1}\circ G_1\circ F_1\\
 &=\Bigl(X+t(t+1)X^2-t^5Y^2-t^3(t+1)^6X^4-2t^3(t+1)XY-2t^2(t+1)^4X^3\\
 &\hspace{1.5cm}-2t^3(t+1)^4X^2Y, Y-t^3(t+1)Y^2-t(t+1)^7X^4\\
 &\hspace{1.5cm}-2t(t+1)^2XY-2(t+1)^5X^3-2t(t+1)^5X^2Y\Bigr)
\end{align*}
Using the algorithrm \ref{algorithrmtotesttameness} we can see that $F\notin \tame_2(R)$.
\end{example}
\section{Structure Of Length Four Automorphisms}
\begin{lemma}\label{C(bX)lemma} Let R be a U.F.D, K its fraction field, and $A(X),\ B(X) \in R[X],\ b \in R$
be such that $A(0)=B(0)=0$, gcd $(B,b)=1$. Then $A({B(X)\over
\displaystyle b}) \in R[X] \Rightarrow A(X)=C(bX)$.
\end{lemma}
  \begin{proof}
\begin{align*} A({B(X)\over b}) =&A(0+{B(X)\over b}) \cr
=&\sum_{i=0}^n A^{(i)}(0){B(X)^i\over i!\ b^i} =\sum_{i=1}^n
A^{(i)}(0){B(X)^i\over i!\ b^i}\cr =&B(X)\sum_{i=1}^n
A^{(i)}(0){B(X)^{i-1}\over i!\ b^i} \in R[X]\cr
\end{align*}
$\Rightarrow B(X)\sum_{i=0}^n A^{(i)}(0)b^{n-i}{\displaystyle B(X)^{i-1}\over i!}\equiv 0\mod b^n$. Since gcd$(B,b)=1$,
we get\\
\begin{equation*}\tag{*}\sum_{i=1}^n A^{(i)}(0)b^{n-i}{B(X)^{i-1}\over \displaystyle i!}\equiv 0\mod b^n
\end{equation*}
 Putting X= 0 in (*) gives us $A^\prime(0)\equiv 0 mod b$.
  i.e Coefficient of X in A(X) is divisible by b. So (*) becomes
  $$\sum_{i=2}^n A^{(i)}(0)b^{n-i}{B(X)^{i-1}\over i!}=B(X)(\sum_{i=2}^n A^{(i)}(0)b^{n-i}{B(X)^{i-1}\over i!})
  \equiv 0\mod b^n)$$
  Again gcd$(B,b)=1$ gives us
\begin{equation*}\tag{**}
\sum_{i=2}^n A^{(i)}(0)b^{n-i}{B(X)^{i-1}\over i!}\equiv 0\mod b^n
\end{equation*}
Putting X=0 in (**) we get, ${\displaystyle A^{\prime
\prime}(0)\over {2!}}\equiv 0 (\text{ mod }\ b^2)$. i.e Coefficient
of $X^2$ in A(X) is divisible by $b^2$. Proceeding like this one
gets that for all $k \ge 1$ coefficient of $X^k$ in A(X) is
divisible by $b^k$.
\end{proof}

Now we prove a lemma about the structure of  automorphisms in $\text{L}^{(4)}(R)$ where $R$ is a U.F.D.
  Let $F =G_2\circ F_2 \circ G_1 \circ F_1 \in \text{L}^{(4)}(R),\ F_i=(X, Y+{\displaystyle A_i(X)\over \displaystyle a_i}),
  \ G_i=(X+{\displaystyle B_i(Y)\over {\displaystyle b_i}},Y),\ A_i(X),\ B_i(X) \in R[X], a_i,\ b_i \in R$ and
  $A_i(0)=B_i(0)=0$ for $i=1,2,\ gcd(A_i,a_i),\ gcd(B_i,b_i)=1$.\\
  \begin{lemma} We use the same notations as above. Then
  $A_2(X)=C(b_1X)\text{ and }B_1(Y)= D(a_2Y)\text{ for some }C(X),\ D(X)\in
  R[X]$ and $ gcd(a_2,b_1)=1$.
  \end{lemma}
  \begin{proof}
\begin{align}
F &=G_2\circ F_2 \circ G_1 \circ F_1 \cr
     &=\Biggl(X+{\displaystyle B_1\bigl( Y+{\displaystyle \displaystyle A_1(X)\over a_1}\bigr) \over \displaystyle b_1}+{\displaystyle B_2\Bigl(Y+{\displaystyle A_1(X)\over \displaystyle a_1}+{A_2 \bigl(X+{\displaystyle B_1(Y+{\displaystyle \displaystyle A_1(X)\over \displaystyle a_1}) \over \displaystyle b_1} \bigr)\over \displaystyle a_2}\Bigr)\over \displaystyle b_2}, \cr
\hspace{.3cm} & \ \ \ \ \ \ \ Y+{\displaystyle A_1(X)\over
a_1}+A_2\Bigl(X+ {\displaystyle B_1\bigl( Y+{\displaystyle
\displaystyle A_1(X)\over a_1}\bigr) \over \displaystyle
b_1}\Bigr)\Biggr)  \end{align}
 Putting X=0 in the second coordinate
of $F$ we get $A_2({\displaystyle B_1(Y)\over \displaystyle b_1}) \in R[X]$.
 Then applying Lemma \ref{C(bX)lemma}, we get $A_2(X)=C(b_1X)$. Similarly putting  $Y=0$ in the first coordinate of
  $F^{-1}$ we get $B_1(Y)=D(a_2Y)$.
  Since
$A_2(X)=C(b_1X)$, we know that gcd$(C(b_1X),a_2)=1
\Rightarrow$ gcd$(a_2,b_1)=1$
  \end{proof}

 \subsection{Proof of the Main Theorem}
Following useful lemma was proved in \cite{TameAndWild} when $R=k[t]$ and was proved when $R$ is a UFD in \cite{Kuttykrishnanpaper1}.
\begin{lemma}\label{length3structurelemma3} Let $F=F_1^{-1}\circ G \circ F_1  \in \text{L}^{(3)}(R)$ where $F_1=(X,Y+\frac{\displaystyle A_1(X)}{\displaystyle a}),\ G=(X+g(Y),Y),\ A_1(X) \in R[X],\ g(Y) \in K[Y],\ a \in R$. Then $g(Y)=D(aY) \text{ for } D(Y) \in R[Y] \text{ and } a \mid D(Y)$.
\end{lemma}
 Let $F_1=(X,Y+\displaystyle \frac {A(X)}{a}), G_1=(X+\displaystyle
\frac {B(Y)}{b}),Y)$ and $F =G_1^{-1}\circ F_1^{-1} \circ G_1 \circ F_1$
where $A(X)\in R[X], B(Y) \in R[Y], A(0)=B(0)=0$.  We may assume
that $gcd(A(X,a))= gcd(B(Y),b)=1$. By Lemma 4 we also know that
$gcd(a,b)=1,\: A(X)=C(bX),\: B(Y)=D(aY)$ for some $C(X)\in R[X],\
D(Y)\in R[Y],$ with $C(0)=D(0)=0$.

\begin{align}\label{lengthfourF}
F=\Bigl( X+ &\frac {D(aY+C(bX))-D(aY+C(bX)-C(bX+D(aY+C(bX))))}{b},\\
  & \ \ Y+ \frac
{C(bX)-C(bX+D(aY+C(bX)))}{a} \Bigr ). \end{align}
\begin{clm}$a\mid D(Y)$
\end{clm}
Let $S=\{1,b,b^2,\dots\}$ and $R_b=S^{-1}R$. Clearly $G_1 \in \aut_2(R_b)$ and so $ G_1\circ F \in \aut_2(R_b)$.
\begin{align*}
G_1\circ F&= G_1\circ G_1^{-1} \circ F_1^{-1} \circ G_1 \circ F_1\\
          &= F_1^{-1}\circ G_1 \circ F_1 \\
          &= (X+\frac{D(aY+C(bX))}{b},Y+ \frac{C(bX)-C(bX+D(aY+C(bX)))}{a} \Bigr ) \in \text{L}^{(3)}(R_b)
\end{align*}
So by Lemma \ref{length3structurelemma3} we have that $a\mid \frac{\displaystyle D(Y)}{\displaystyle b}$ in $R_b$. Since $gcd(a,b)=1$, this implies that $a\mid D(Y)$. Hence the claim.
.

Notice that by the claim $a$ divides all the terms in the first coordinate of $F$ in \ref{lengthfourF} except $X$. So we have
\begin{align*}
F &=(X+aP(X,Y), Y+Q(X,Y))\ \text{ where }\\
P(X,Y)&=\frac {\displaystyle
D(aY+C(bX))-D(aY+C(bX)-C(bX+D(aY+C(bX))))}{\displaystyle
 ab}\ \text{ and }\\
  Q(X,Y)&={\displaystyle \frac {C(bX)-C(bX+D(aY+C(bX)))} {a}}.
  \end{align*}
\

Let $E=(X+aW,Y,W)$ and $L=(X,Y-Q(X,0),W-P(X,0))$ then
\begin{align*}
F & \stameequiv (X+aP(X,Y), Y+Q(X,Y),W) \\
& \tameequiv F_1= (X+aP(X,Y),Y+Q(X,Y),W+P(X,Y))\\
& \tameequiv E \circ F_1 \circ E^{-1}= (X,Y+ Q(X+aW,Y),
W+P(X+aW,Y)). \\
& \tameequiv F^1= L \circ E \circ F_1 \circ
 E^{-1}\\
&=(X,Y+Q(X+aW,Y)-Q(X,0),W+P(X+aW,Y)-P(X,0))\\
=\Biggl(X,Y+&{\displaystyle \frac {C(bX+abW)-C(bX+abW+D(aY+C(bX+abW)))} {a}}\\
&\hspace{5cm}-\frac{C(bX)+C(bX+D(C(bX)))}{a},\\
W+&\frac {\displaystyle D(aY+C(bX+abW))-D(C(bX))}{\displaystyle ab}\\
-&\frac{D(aY+C(bX+abW)-C(bX+abW+D(aY+C(bX+abW))))}{\displaystyle ab}\\
 &\hspace{5.5cm} +\frac{D(C(bX)-C(bX+D(C(bX))))}{ab}\Biggr)
\end{align*}
We can compute that $F^1=G_2^1 \circ F_{2}^1 \circ G_1^1 \circ F_{1}^1$
where
\begin{align*}
F_1^1&=(X,Y+{{\displaystyle C(b(X+aW))-C(bX)\over a}},W),\\
G_1^1&=(X,Y, W+{{\displaystyle D(aY+C(bX))-D(C(bX))\over ab}}),\\
F_2^1&=(X,Y-{{\displaystyle C(b(X+aW)+D(C(bX)))+C(D(C(bX)))\over a}},W)\\
\text{and }G_2^1&=\Bigl(X,Y,W-\frac{\displaystyle D(bY+C(bX)-C(bX+D(C(bX))))}{ab}\\
&\hspace{4cm} + \frac{\displaystyle D(C(bX)-C(bX+D(C(bX))))}{ab}\Bigr)
\end{align*}
 Since $F_{1}^1
\in \El_2(R[X])$, $F^1$ is tamely equivalent to $F^1 \circ
(F_{1}^1)^{-1} \in \text{L}^{(3)}(R[X])$ and hence stably tame by Theorem
\ref{length3isstablytame1}. Thus we get that $F$ is stably tame.
\subsection{An intriguing Example}
Using the notations above we give an example of a length four automorphism which is not a commutator. Further, in this example, $a \neq  b$. However, this automorphism is stably tame!
\begin{example}\label{lengthfourexample}
  Let $R$ be a domain and $t\in R\backslash \{ 0\}.$
  \begin{align*}
  F_1 &= (X,Y+\frac{X^2}{t})\\
  F_2 &= (X,Y+(t-1)X)\\
  G_1 &= (X+(t+1)Y,Y)\\
  G_2 &= (X-\frac{Y^2}{t},Y) \text{ and } \\
  F &= G_2\circ F_2 \circ G_1 \circ F_1 \\
    &=(X+(t+1)Y+3X^2-t^3Y^2-tX^2-tX^4-2t^2XY\\
    &\hspace{1cm}+2tXY-2t^2X^2Y-2tX^3+2X^3,\\
    &\hspace{5.5cm}t^2Y+(t-1)X+tX^2)
  \end{align*}
 Then $F\in \text{L}^{(4)}(R)$. Using the Algorithrm \ref{algorithrmtotesttameness} we get that $F \notin \tame_2(R)$.
  \end{example}
 \begin{align*}\text{Let }P(X,Y)&= (t+1)Y+3X^2-t^3Y^2-tX^2-tX^4-2t^2XY+2tXY\\
                                           &\hspace{5.5cm}-2t^2X^2Y-2tX^3+2X^3 \\
                         Q(X,Y)&= t^2Y+(t-1)X+tX^2 \\
\text{ and } \widetilde{Q}(X,Y)&=X+tY+X^2.
\end{align*}
We extend $F$ to $\widetilde{F}=(F,Z)\in \aut_3(R)$ and define the following elementary automorphisms of $R[X,Y,Z]$.
\begin{align*}
\tau =&(X,Y, W+\widetilde{Q}(X))\\
\eta =&(X,Y-tZ,Z)\\
\phi =&(X-tZ,Y,Z)\\
\text{Also, let }\pi=&(-Y,X,Z)\\
\text{Then }\pi \circ \eta \circ \widetilde{F} \circ \tau =&(X+tZ,X+P(X,Y),Z+\widetilde{Q}(X,Y))\text{ and} \\
&\hspace{-3.5cm}\widetilde{F^1}=\pi \circ \eta \circ F^1 \circ \tau \circ \phi\\
  &\hspace{-3cm}=(X,X-tZ+P(X-tZ,Y),Z+\widetilde{Q}(X-tZ,Y))\\
 &\hspace{-3cm}=\Bigl(X, X-tZ+ (t+1)Y+t^3Y^2+3(X-tZ)^2-t(X-tZ)^2\\
 & \hspace{-2cm} -t(X-tZ)^4-2t^2(X-tZ)Y+2t(X-tZ)Y-2t^2(X-tZ)^2Y\\
 & \hspace{-1cm} -2t(X-tZ)^3+2(X-tZ)^3, Z+ (X-tZ)+tY+(X-tZ)^2 \Bigr)\\
&\hspace{-3cm}=\Bigl(X,Y+tY+ X-tZ+ t^3Y^2+3(X-tZ)^2-t(X-tZ)^2-t(X-tZ)^4\\
&\hspace{-1cm}-2t^2(X-tZ)Y+2t(X-tZ)Y-2t^2(X-tZ)^2Y\\
& \hspace{-1cm}-2t(X-tZ)^3+2(X-tZ)^3, Z+ (X-tZ)+tY+(X-tZ)^2\Bigr)\\
 \end{align*}
\begin{align*}
\text{Then }\widetilde{F^1}=&(X, Y+P_1(X,Y,Z), Z+ Q_1(X,Y,Z))\text{ where } \\
P_1(X,Y,Z)=&P(X-tZ,Y)-Y+X \text{ and }Q_1(X,Y,Z)=\widetilde{Q}(X-tZ,Y)\\
\end{align*}
Notice that $\widetilde{F}\stameequiv \widetilde{F^1}$. Let $\Theta=(X,Y-P_1(X,0,0), Z-Q_1(X,0,0))$.\\
 Then $\widetilde{F^1}\tameequiv \Theta \circ \widetilde{F^1}$.
Clearly the following automorphisms are in $\BA_3(R)$.
\begin{align*}
\widetilde{F_1}&=(X,Y+\frac{(X+tZ)+(X+tZ)^2-X-X^2}{t},Z)\\
\widetilde{G_1}&=(X,Y,Z+tY)
\end{align*}
Then we have that $\Theta \circ \widetilde{F^1}=\widetilde{F_1}^{-1}\circ \widetilde{ G_1} \circ  \widetilde{F_1}$. So $\tau_1 \circ \widetilde{F^1} \in \tame_3(R)$. Hence $F$ is stably tame.
\bibliographystyle{alpha}
\bibliography{stablytamepaper}

\end{document}